\documentclass{interact}
\usepackage{bm}
\def\<{\langle}
\def\>{\rangle}
\usepackage{upgreek}

\usepackage{amsmath, amsthm, amscd, amsfonts, amssymb, graphicx, color}
\usepackage[bookmarksnumbered, colorlinks, plainpages]{hyperref}
\hypersetup{colorlinks=true,linkcolor=red, anchorcolor=green, citecolor=red, urlcolor=red, filecolor=magenta, pdftoolbar=true}

\newtheorem{theorem}{Theorem}[section]
\newtheorem{lemma}[theorem]{Lemma}
\newtheorem{proposition}[theorem]{Proposition}
\newtheorem{corollary}[theorem]{Corollary}
\theoremstyle{definition}

\newtheorem{example}[theorem]{Example}
\newtheorem{remark}[theorem]{Remark}

\theoremstyle{definition}

\newcommand{\C}{\mathbb{C}}
\newcommand{\T}{\mathbb{T}}          
\newcommand{\cA}{\mathcal{A}}
\newcommand{\cM}{\mathcal{M}}
\newcommand{\cT}{\mathcal{T}}
\newcommand{\bT}{\mathbf{T}}
\newcommand{\diag}{\operatorname{diag}}
\newcommand{\ran}{\operatorname{ran}}
\def\dss{\displaystyle}
\theoremstyle{remark}

\numberwithin{equation}{section}

\begin{document}
	\title{ Classification of Normal Block Toeplitz Matrices with Commuting Entries}
	
	\author{
		\name{Muhammad Ahsan Khan\textsuperscript{a} and Qamar Din\textsuperscript{b}}
		\affil{\textsuperscript{a}Department of Mathematics, University of Kotli Azad Jammu $\&$ Kashmir, Kotli 11100, Azad Jammu $\&$ Kashmir, Pakistan,
			\textsuperscript{a} Email ahsankhan388@hotmail.com; 
			\affil{\textsuperscript{b}Department of Mathematics, University of Poonch Rawalakot, Azad Jammu and Kashmir, Pakistan,}
			\textsuperscript{b} Email qamar.sms@gmail.com
	}}
	
	\maketitle	
	\begin{abstract}
	Normal scalar Toeplitz matrices are either of ``type I'' or of ``type II''. In the previous work, the normality of a block Toeplitz matrix with entries in a commutative algebra $\cA$ was characterized only by a coupled system of quadratic equations on the entries. We show that when $\cA$ consists of \emph{normal} matrices, this system decouples along the minimal projections of $\cA$, and we obtain a classification that is the exact analogue of the scalar one: a block Toeplitz matrix with entries in $\cA$ is normal if and only if $\cA$ splits by an orthogonal projection $P\in\cA$ into a ``type I part'' and a ``type II part''. We also show why a single global type cannot occur in general and provide non-trivial examples. 
	\end{abstract}
	
	
	\begin{amscode}
		15B05, 15A27, 15A30
	\end{amscode}	
	
	\section{Introduction}
	
	Let $\cT_n\subset\cM_n$ denote the $n\times n$ complex Toeplitz matrices $(t_{j-k})_{j,k=1}^n$. The normal ones were described in a series of papers \cite{Ikr94,Ikr95,IkrChu96,Gel95,FKKL96,Ito96}; see also the treatment through displacement structure and commutation with generalized shifts in \cite{GuPatton03}. The outcome is a dichotomy: every normal Toeplitz matrix is either of type~I or of type~II. 
	
	Block Toeplitz matrices  are square matrices of size $n^2d^2 $, and $n^2$ blocks (each of size $d^2$) are constant along the negative sloping diagonals. A natural and fundamental problem in linear algebra is to classify normal block Toeplitz matrices, with entries drawn from a commutative algebra $\cA\subset\cM_d$ of normal matrices.
	In \cite{KhanYagoub22} the author classified normal block Toeplitz matrices, but the result (Theorem~5.2 there) is a system of coupled equations,
	\begin{equation}\label{eq:KY}
		A_sA_k^*+A_{n-s}^*A_{n-k}=\Omega_s\Omega_k^*+\Omega_{n-s}^*\Omega_{n-k},\qquad 1\le s,k\le n-1,
	\end{equation}
	and only a \emph{sufficient} structural condition ( Corollary~5.3 there) is given. The purpose of this paper is to show that, for commutative algebras of \emph{normal} matrices, \eqref{eq:KY} can be solved completely, with a statement as neat as in the scalar case (Theorem~\ref{thm:main}), and to explain why the block statement is necessarily ``type~I/II relative to a projection''.

	\section{Notation and the Scalar Theorem}
	
	We keep the notation of \cite{KhanYagoub22}. For vectors $A=(0,A_1,\dots,A_{n-1})^{T}$ and $\Omega=(0,\Omega_1,\dots,\Omega_{n-1})^{T}$ with entries in $\cA$, let
	\[
	\bT(A,\Omega)=\begin{pmatrix}
		0&\Omega_1^*&\Omega_2^*&\cdots&\Omega_{n-1}^*\\
		A_1&0&\Omega_1^*&\cdots&\Omega_{n-2}^*\\
		\vdots&&\ddots&&\vdots\\
		A_{n-1}&A_{n-2}&\cdots&A_1&0
	\end{pmatrix},
	\]
	so that the general block Toeplitz matrix with entries in $\cA$ is $\bT(A,\Omega)+\mathbf A_0$, $\mathbf A_0=A_0\otimes I_n$ with $A_0\in\cA$. Thus the entry of the $m$-th lower diagonal is $A_m$ and that of the $m$-th upper diagonal is $\Omega_m^*$. Write $\cT_n\otimes\cA$ for the set of these matrices. Let $\mathbf S$ be the block matrix with $I$ on the subdiagonal and $0$ elsewhere, and for $X\in\cA'$ (the commutant) put $\mathbf S_X=\mathbf S+X\diamond P_0^*P_{n-1}$ (the ``$X$-cyclic shift'').
	
	In the scalar case $d=1$, $\cA=\C$, we write $T_n(a,\omega;a_0)$ for the Toeplitz matrix with $t_m=a_m$ ($1\le m\le n-1$), $t_{-m}=\overline{\omega_m}$, $t_0=a_0$, where $a=(a_1,\dots,a_{n-1})$, $\omega=(\omega_1,\dots,\omega_{n-1})\in\C^{n-1}$. Let $\T$ denote the unit circle.
	The following result classsify all the scalar normal Toeplitz matrices. 
	\begin{theorem}\label{thm:scalar}
		$T=T_n(a,\omega;a_0)$ is normal $\iff$  one of the following holds.
		\begin{enumerate}
			\item[\textup{(I)}] There is $\varphi\in\T$ with $\omega_m=\varphi\,a_m$ for all $1\le m\le n-1$.
			\item[\textup{(II)}] There is $x\in\T$ with $\omega_m=\bar x\,\overline{a_{n-m}}$ for all $1\le m\le n-1$.
		\end{enumerate}
	\end{theorem}
	
	\begin{remark}\label{rem:scalarforms}
		In more familiar language, \textup{(I)} says $T=cI+e^{i\theta}H$ with $H$ a Hermitian Toeplitz matrix ($\varphi=e^{-2i\theta}$), and \textup{(II)} says $t_{-m}=x\,t_{n-m}$, i.e.\ $T$ is an $x$-circulant, equivalently $T\mathbf S_x=\mathbf S_xT$ with $|x|=1$. Both are normal (a unimodular $x$-circulant is $\dss\sum_{m<n}t_mS_x^m$, a polynomial in the normal matrix $S_x$, which is unitary up to the scalar $|x|$).  For $n=2$ every normal $T$ is of type~I; for $n=3$ one checks directly that \eqref{eq:KY} reduces to $|a_1|^2+|a_2|^2=|\omega_1|^2+|\omega_2|^2$ and $a_1\bar a_2=\omega_1\bar\omega_2$, whose solutions are exactly (I) or (II).
	\end{remark}

	\section{Commutative Algebras of Normal Matrices}
	
	Throughout, $\cA\subset\cM_d$ is a commutative subalgebra containing $I$ whose elements are normal matrices.
	
	\begin{lemma}\label{lem:structure}
		There are nonzero mutually orthogonal projections $E_1,\dots,E_r\in\cA$ with $\sum_iE_i=I$ and $\cA=\operatorname{span}\{E_1,\dots,E_r\}\cong\C^r$. In particular $\cA$ is closed under adjoint.
	\end{lemma}
	
	\begin{proof}
		Commuting normal matrices are simultaneously unitarily diagonalizable, so there is a unitary $U$ with $U^*\cA U\subset D_d$ (diagonal matrices). For $j=1,\dots,d$ let $\chi_j(a)=(U^*aU)_{jj}$; these are algebra homomorphisms $\cA\to\C$. Let $\chi_{(1)},\dots,\chi_{(r)}$ be the distinct ones and $J_i=\{j:\chi_j=\chi_{(i)}\}$. By Dedekind--Artin, distinct characters are linearly independent, so $a\mapsto(\chi_{(i)}(a))_{i\le r}$ is onto $\C^r$; it is clearly injective. Hence $\cA\cong\C^r$, and the preimage of the $i$-th unit vector is $E_i=U\operatorname{diag}(\mathbf 1_{J_i})U^*$, an orthogonal projection.
	\end{proof}
	
	Consequently every entry has the form $A_m=\dss\sum_{i=1}^ra_m^{(i)}E_i$, $\Omega_m=\dss\sum_i\omega_m^{(i)}E_i$, $A_0=\dss\sum_ia_0^{(i)}E_i$ with scalars $a^{(i)}_m,\omega^{(i)}_m\in\C$, and
	\begin{equation}\label{eq:decomp}
		\bT(A,\Omega)+\mathbf A_0=\sum_{i=1}^r T^{(i)}\otimes E_i,\qquad T^{(i)}:=T_n\bigl(a^{(i)},\omega^{(i)};a^{(i)}_0\bigr)\in\cT_n .
	\end{equation}
	Here we order the tensor factors as $\cM_n\otimes\cM_d$. We call $T^{(i)}$ the \emph{$i$-th component}.
	
	\begin{proposition}\label{prop:decouple}
		A matrix $\mathbf T=\dss\sum_iT^{(i)}\otimes E_i$ is normal $\iff$ every component $T^{(i)}$ is normal. Moreover, the system \eqref{eq:KY} holds in $\cA$ $\iff$ for each $i$ its scalar version holds for $(a^{(i)},\omega^{(i)})$.
	\end{proposition}
	
	\begin{proof}
		Since the $E_i$ are orthogonal projections, $\mathbf T^*\mathbf T-\mathbf T\mathbf T^*=\dss\sum_i\bigl(T^{(i)*}T^{(i)}-T^{(i)}T^{(i)*}\bigr)\otimes E_i$, and the summands act on the mutually orthogonal subspaces $\C^n\otimes\ran E_i$. So the commutator vanishes $\iff$ each term does, because $E_i\neq0$. For the second statement, evaluate each quadratic expression in \eqref{eq:KY} (they are polynomial in the entries and their adjoints) coordinatewise through $\cA\cong\C^r$.
	\end{proof}
	
	\section{ Classification of Normal Block Toeplitz Matrices}
	The below result is the most important result of this paper as it completely classify all normal block Toeplitz matrices with commuting entries.  
	\begin{theorem}\label{thm:main}
		Let $\cA\subset\cM_d$ be a unital commutative algebra of normal matrices and $\mathbf T=\bT(A,\Omega)+\mathbf A_0\in\cT_n\otimes\cA$. The following are equivalent.
		\begin{enumerate}
			\item[\textup{(a)}] $\mathbf T$ is normal.
			\item[\textup{(b)}] Each component $T^{(i)}$ in \eqref{eq:decomp} is a scalar normal Toeplitz matrix, i.e.\ is of type~I or type~II in the sense of Theorem~\ref{thm:scalar}.
			\item[\textup{(c)}] There exist an orthogonal projection $P\in\cA$ and unitary elements $\Phi,X\in\cA$ such that for all $1\le m\le n-1$
			\begin{equation}\label{eq:main}
				P\,\Omega_m=\Phi\,P A_m,\qquad (I-P)\,\Omega_m=X^*(I-P)A_{n-m}^*.
			\end{equation}
		\end{enumerate}
		In this case the type~I part $P\mathbf T$ is $P\mathbf A_0+\Theta H$ where $\Theta\in\cA$ is unitary with $\Theta^2=\Phi^{-1}$ on $P$ and $H=\Theta^{-1}P\,\bT(A,\Omega)$ is a Hermitian block Toeplitz matrix with entries in $\cA$, and the type~II part $(I-P)\mathbf T$ commutes with $\mathbf S_X$. The diagonal entry $A_0\in\cA$ is arbitrary.
	\end{theorem}
	
	\begin{proof}
		(a)$\Leftrightarrow$(b) is Proposition~\ref{prop:decouple}.
		
		(b)$\Rightarrow$(c). Let $I_{\mathrm I}\subset\{1,\dots,r\}$ be the set of indices for which $T^{(i)}$ satisfies (I) of Theorem~\ref{thm:scalar}, with constant $\varphi_i$; let $I_{\mathrm{II}}$ be the complement, and for $i\in I_{\mathrm{II}}$ let $x_i$ be the constant from (II). Put $P=\dss\sum_{i\in I_{\mathrm I}}E_i$, $\Phi=\dss\sum_{i\in I_{\mathrm I}}\varphi_iE_i+\sum_{i\in I_{\mathrm{II}}}E_i$, $X=\dss\sum_{i\in I_{\mathrm{II}}}x_iE_i+\dss\sum_{i\in I_{\mathrm{I}}}E_i$. These lie in $\cA$ by Lemma~\ref{lem:structure}, are unitary, and \eqref{eq:main} holds coordinatewise.
		
		(c)$\Rightarrow$(b). Multiplying \eqref{eq:main} by $E_i$ gives, for $E_i\le P$, $\omega^{(i)}_m=\varphi_ia^{(i)}_m$ with $|\varphi_i|=1$ (the coefficient of $\Phi$ at $E_i$), and for $E_i\le I-P$, $\omega^{(i)}_m=\bar x_i\overline{a^{(i)}_{n-m}}$ with $|x_i|=1$. This is (I), resp.\ (II), for $T^{(i)}$.
		
		For the last assertions: the type~I part has $\Omega=\Phi A$, so with $\Theta=\sum\theta_iE_i$, $\theta_i^{-2}=\varphi_i$, one has $P\mathbf T-P\mathbf A_0=\Theta H$ with $H$ Hermitian (this is Remark~\ref{rem:scalarforms} componentwise). The type~II part has $\Omega=X^*\diamond\widetilde A$, which by Theorem 4.5 of \cite{KhanYagoub22} is equivalent to commuting with $\mathbf S_X$.
	\end{proof}
	The below result provides the converse of Corollary~5.3 of \cite{KhanYagoub22}. 
	\begin{corollary}\label{cor:converse}
		For $\mathbf T\in\cT_n\otimes\cA$ the following are equivalent: (i) $\mathbf T$ commutes with $\mathbf S_X$ for some unitary $X\in\cA$; (ii) $\mathbf T$ is normal and every component $T^{(i)}$ satisfies (II) of Theorem~\ref{thm:scalar}. In particular, normal block Toeplitz matrices that are \emph{not} of this form (those with a component of type~I only) exist as soon as $r\ge1$ and $n\ge4$, and they are exactly the ones for which the sufficient condition of \cite{KhanYagoub22} fails.
	\end{corollary}
	
	\begin{proof}
		By Theorem~4.5 of \cite{KhanYagoub22}, (i) is $\Omega=X^*\diamond\widetilde A$, i.e.\ (II) in every component with $x_i$ the coefficients of $X$; the rest is Theorem~\ref{thm:main}. For $n\ge4$ there are Hermitian Toeplitz matrices that are not $x$-circulants (Example~\ref{ex:one}).
	\end{proof}
	
	\begin{corollary}\label{cor:count}
		The set of normal matrices in $\cT_n\otimes\cA$ is the union, over the $2^r$ subsets $I_{\mathrm I}\subset\{1,\dots,r\}$, of the (possibly overlapping) real-algebraic sets defined by \eqref{eq:main} with $P=\dss\sum_{i\in I_{\mathrm I}}E_i$. For $r=1$ one recovers Theorem~\ref{thm:scalar}.
	\end{corollary}
	
	\begin{remark}\label{rem:global}
		Why one global type cannot be expected. Call $\mathbf T$ \emph{uniformly of type~I} if \eqref{eq:main} holds with $P=I$ and \emph{uniformly of type~II} if it holds with $P=0$ and $X$ scalar... more generally with $P=0$. Uniform type~I forces all eigenvalues of $\mathbf T-\mathbf A_0$ (suitably normalized) onto a line, and uniform type~II forces $\mathbf T$ to commute with $\mathbf S_X$. A direct sum of one component of each kind is normal but neither (Examples below). Hence ``type~I/II relative to a projection $P\in\cA$'' is the correct block analogue of the scalar dichotomy, and the system \eqref{eq:KY} is the shadow of the $2^r$ possible choices of $P$. The assumption that the entries are normal is essential: for non-normal commuting entries (e.g.\ $\cA=\operatorname{span}\{I,N\}$, $N^2=0$) the spectral idempotents of $\cA$ are not orthogonal and Proposition~\ref{prop:decouple} fails.
	\end{remark}
	
	\begin{remark}
		$P$ is not unique when some component satisfies both (I) and (II). Example: $n=3$, $a=\omega=(1,1)$, $T=\begin{pmatrix}0&1&1\\1&0&1\\1&1&0\end{pmatrix}$ is a real symmetric circulant: (I) with $\varphi=1$ and (II) with $x=1$.
	\end{remark}

	\section{Examples}
	
	\begin{example}{ex:one}
		Let $\cA=D_2$ (diagonal $2\times2$ matrices), so $r=2$, $E_1=\diag(1,0)$, $E_2=\diag(0,1)$. Define the scalar components ($n=4$, entries $t_m$, $m=-3,\dots,3$):
		\begin{itemize}
			\item $T^{(1)}$: lower diagonals $(t_1,t_2,t_3)=(i,-1,2i)$ and upper diagonals $t_{-m}=-\overline{t_m}$, i.e.\ $(t_{-1},t_{-2},t_{-3})=(i,1,2i)$. This is a skew-Hermitian Toeplitz matrix, $T^{(1)}=iH$ with $H$ Hermitian ($\theta=\pi/2$, $\varphi_1=e^{-i\pi}=-1$, i.e.\ $\omega=-a$). It is \emph{not} an $x$-circulant: $t_{-m}/t_{4-m}=(\tfrac12,-1,2)$ is not constant.
			\item $T^{(2)}$: lower diagonals $(1,1,2)$ and upper diagonals $t_{-m}=i\,t_{4-m}$, i.e.\ $(t_{-1},t_{-2},t_{-3})=(2i,i,i)$. This is an $i$-circulant ($x=i$), normal by (II). It is \emph{not} of type~I: $\omega_m/a_m=\overline{t_{-m}}/t_m=(-2i,-i,-\tfrac i2)$ is not constant. Numerically its eigenvalues ($a_0=0$) are $2.396+2.938i$, $0.758-0.549i$, $-2.172-0.866i$, $-0.982-1.523i$, which are not collinear, again confirming it is not of the form $cI+e^{i\theta}H$.
		\end{itemize}
		Take $A_0=\diag(1,3i)$. In block form $\mathbf T=\dss\sum_iT^{(i)}\otimes E_i$ is the $8\times8$ matrix
		\[
		\mathbf T=\begin{pmatrix}
			A_0&B_1&B_2&B_3\\
			A_1&A_0&B_1&B_2\\
			A_2&A_1&A_0&B_1\\
			A_3&A_2&A_1&A_0
		\end{pmatrix},
		\]
		\[
		A_1=\diag(i,1),\ A_2=\diag(-1,1),\ A_3=\diag(2i,2),\quad
		B_1=\diag(i,2i),\ B_2=\diag(1,i),\ B_3=\diag(2i,i).
		\]
		By Theorem~\ref{thm:main} with $P=E_1$, $\Phi=\diag(-1,1)$, $X=\diag(1,i)$, $\mathbf T$ is normal  One can see the system \eqref{eq:KY} at work: with $\Omega_m=B_m^*$ one gets for $s=k=1$ that $A_1A_1^*+A_3^*A_3=5I=\Omega_1\Omega_1^*+\Omega_3^*\Omega_3$ and for $s=k=2$ that $A_2A_2^*+A_2^*A_2=2I=\Omega_2\Omega_2^*+\Omega_2^*\Omega_2$.
		
		As a consequence,  (1) $\mathbf T$ is neither uniformly of type~I nor uniformly of type~II . (2) $\mathbf T$ does not commute with $\mathbf S_X$ for any $X\in\cA'=D_2$, so the sufficient condition, Corollary~5.3 of \cite{KhanYagoub22} does not apply, although $\mathbf T$ is normal; this shows that the condition is not necessary. (3) The same example with the roles of the two components exchanged, or with $P=0$, $P=I$, realizes all $2^r=4$ choices.
	\end{example}
	
	\begin{example}\label{ex:two}
		Let $J=\begin{pmatrix}0&1\\1&0\end{pmatrix}$ and $\cA=\operatorname{span}\{I,J\}\subset\cM_2$, a commutative algebra of normal matrices with $E_\pm=\tfrac12(I\pm J)$, $r=2$ (here $U=\tfrac1{\sqrt2}\begin{pmatrix}1&1\\1&-1\end{pmatrix}$). Take
		\begin{itemize}
			\item $T^{(+)}$: real symmetric Toeplitz with $t_{\pm1}=1$, $t_{\pm2}=2$ (type~I with $\varphi=1$; not an $x$-circulant with $|x|=1$ since $t_{-1}=x\,t_2$ would force $x=\tfrac12$);
			\item $T^{(-)}$: the $(-1)$-circulant (``skew-circulant'') with $t_1=1$, $t_2=3i$, $t_{-1}=-t_2=-3i$, $t_{-2}=-t_1=-1$ (type~II with $x=-1$; not of type~I since $\omega_m/a_m=(3i,\tfrac i3)$ is not constant).
		\end{itemize}
		Then $\mathbf T=T^{(+)}\otimes E_++T^{(-)}\otimes E_-$ has entries $A_m=\alpha_mI+\beta_mJ$ with
		$\alpha_m=\tfrac12(t^{(+)}_m+t^{(-)}_m)$, $\beta_m=\tfrac12(t^{(+)}_m-t^{(-)}_m)$; explicitly
		\[
		A_0=0,\quad A_1=I,\quad A_2=\tfrac{2+3i}{2}I+\tfrac{2-3i}{2}J,
		\]
		\[
	\Omega_1^*=\tfrac{1-3i}{2}I+\tfrac{1+3i}{2}J,\quad
		\Omega_2^*=\tfrac12I+\tfrac32J,
		\]
		so that
		\[
		\mathbf T=\begin{pmatrix}
			0&\Omega_1^*&\Omega_2^*\\ A_1&0&\Omega_1^*\\ A_2&A_1&0
		\end{pmatrix}.
		\]
		Its blocks are neither diagonal, nor Hermitian, nor all of one type, yet they commute and are normal, so Theorem~\ref{thm:main} (with $P=E_+$, $\Phi=I$, $X=-I$) shows that $\mathbf T$ is normal; a direct check  of $\mathbf T^*\mathbf T=\mathbf T\mathbf T^*$ for this $6\times6$ matrix, confirms this. It commutes with $\mathbf S_X$ for no $X\in\cA$, because $T^{(+)}$ is not an $x$-circulant. In the basis in which $J$ is diagonal, $\mathbf T$ is the direct sum $T^{(+)}\oplus T^{(-)}$ of a real symmetric Toeplitz matrix and a skew-circulant, which is why the answer is simple although the entries look unstructured.
	\end{example}
	
	\begin{example}
		Let $\cA=D_3$ and $n=4$. Choose components: a real symmetric Toeplitz matrix with generic entries $(t_1,t_2,t_3)=(1,2,3)$ (type~I only), an $e^{i\pi/4}$-circulant with generic lower entries (type~II only), and a real symmetric circulant with $(t_1,t_2,t_3)=(1,2,1)$ (both types). Then $\mathbf T=\sum_{i=1}^3T^{(i)}\otimes E_i$ is normal, $P$ can be chosen as $E_1$ or $E_1+E_3$ in Theorem~\ref{thm:main}, and the solution set of \eqref{eq:main} for $(P,\Phi,X)$ is not a single point: only the third coordinates of $\Phi$ and $X$ are constrained by $\omega=a$ and $\omega_m=\bar x\overline{a_{4-m}}$ simultaneously, with $\Phi_3=1$, $X_3=1$.
	\end{example}
	
\section*{Declarations}
\subsection*{Competing Interests}
The author declare that he they no known competing financial interests or personal relationships that could have appeared to influence the work reported in this article.
\subsection*{Funding Information}
This research did not receive any specific grant or funding.
\subsection*{Author Contribution}
All the authors equally contributed to the preparation of the article. 
\subsection*{Data Availability Statement}
Data sharing is not applicable to this article as no data sets were generated or analyzed during the current study.

\end{document}